\documentclass{amsart}

\newtheorem{theorem}{Theorem}[section]
\newtheorem{lemma}[theorem]{Lemma}
\newtheorem{claim}[theorem]{Claim}
\newtheorem{corollary}[theorem]{Corollary}
\newtheorem{proposition}[theorem]{Proposition}
\usepackage{graphicx}

\theoremstyle{definition}

\theoremstyle{remark}
\newtheorem{remark}[theorem]{Remark}
\newtheorem{question}[theorem]{Question}

\numberwithin{equation}{section}

\begin{document}

\title[Knots yielding homeomorphic lens spaces]{Knots yielding homeomorphic lens spaces by Dehn surgery}
\author{Toshio Saito}
\address{Department of Mathematics, 
Graduate School of Humanities and Sciences, Nara Women's University,
Kita-Uoya Nishimachi, Nara 630-8506, Japan}
\email{tsaito@cc.nara-wu.ac.jp}
\thanks{The first author was supported by the 21st Century COE program \lq\lq Towards a New Basic Science; Depth and Synthesis\rq\rq, Osaka University.}
\author{Masakazu Teragaito}
\address{Department of Mathematics and Mathematics Education, Hiroshima University,
1-1-1 Kagamiyama, Higashi-hiroshima 739-8524, Japan}
\email{teragai@hiroshima-u.ac.jp}
\thanks{The second author was partially supported by Japan Society for the Promotion of Science,
Grant-in-Aid for Scientific Research (C), 19540089.}

%    General info

%\date{January 1, 1994 and, in revised form, June 22, 1994.}

%\dedicatory{This paper is dedicated to our authors.}


\begin{abstract}
We show that there exist infinitely many pairs of distinct knots in the $3$-sphere such that
each pair can yield homeomorphic lens spaces by the same Dehn surgery. 
Moreover, each knot of the pair can be chosen to be a torus knot, a satellite knot or a hyperbolic knot,
except that both cannot be satellite knots simultaneously.
This exception is shown to be unavoidable by the classical theory of binary quadratic forms.
\end{abstract}

\maketitle

%%%%%%%%%%%%%%%%%%%%%%%%%%%%%%%%%%%%%%%%%%%%%%%%%%%%%%%%%%%%%%%%%%%%%%%%
\section{Introduction}

For a knot $K$ in the $3$-sphere $S^3$, let $K(m/n)$ denote
the closed oriented $3$-manifold obtained by $m/n$-Dehn surgery on $K$, which
is the union of the knot exterior $E(K)=S^3-\mathrm{int}N(K)$ and a solid torus $V$
in such a way that the meridian of $V$ is attached to a loop on $\partial E(K)$ with slope $m/n$.
In this paper, all $3$-manifolds are oriented, and
two knots in $S^3$ are said to be \textit{equivalent\/} if there is an orientation-preserving
homeomorphism of $S^3$ sending one to the other.

For a fixed slope $m/n$, $m/n$-surgery can be regarded as a map from the set of the equivalence classes of knots to
that of $3$-manifolds.
There are many results on the injectivity of this map.
Lickorish \cite{L} gave two non-equivalent knots on which $(-1)$-surgeries yield
the same homology sphere.
Brakes \cite{Br} showed that for any integer $n\ge 2$, there exist $n$ distinct knots on which
$1$-surgeries yield the same $3$-manifold.
See also \cite{Ka,Liv,Te}.
Finally, Osoinach \cite{O1,O2} showed the existence of $3$-manifolds, in fact, a hyperbolic $3$-manifold and
a toroidal manifold, which can be obtained from infinitely many hyperbolic knots
by $0$-surgery.
By using Osoinach's construction, the second author gave a Seifert fibered manifold over the $2$-sphere with three
exceptional fibers which can be obtained from infinitely many hyperbolic knots
by $4$-surgery \cite{Te2}.
Thus it is natural to ask whether there exists a lens space which can be obtained from infinitely many knots
by the same Dehn surgery or not.
Although we do not know the answer to it yet, we feel it negative through our computer experiment.
In fact, as far as we know, at most two knots can yield homeomorphic lens spaces by the same Dehn surgery.

We should note that Berge's table \cite{B} shows that
there are $32$ lens spaces, among those with fundamental groups of order up to $500$, which admit two knots yielding $S^3$ by Dehn surgery.
This strongly suggests that many lens spaces can be obtained from non-equivalent knots in $S^3$ by the same Dehn surgery.
In this paper, we study whether a pair of non-equivalent knots can yield homeomorphic lens spaces, ignoring orientations,
by the same Dehn surgery.
We should be attentive to this orientation convention.
Let $U$ be the unknot and $K$ a knot in $S^3$.
By using Floer homology for Seiberg-Witten monopoles, it is proved in \cite{KMOS}
that if there exists an orientation-preserving homeomorphism between
$K(m/n)$ and $U(m/n)$ then $K$ is trivial.
In other words, if $K(m/n)$ is homeomorphic to the lens space $L(m,n)$ under an orientation-preserving
homeomorphism, then $K$ is trivial.
Here, the preservation of orientation is important, because $5$-surgery on the right-handed trefoil yields
$L(5,4)=L(5,-1)$. From our point of view, the right-handed trefoil and the unknot yield
homeomorphic lens spaces under the same $5$-surgery.

As a consequence of the cyclic surgery theorem \cite{CGLS}, any non-trivial amphicheiral knot
has no Dehn surgery yielding a lens space, and
the pair of a knot and its mirror image cannot yield homeomorphic lens spaces by the same Dehn surgery.
Also, only torus knots admit non-integral lens space surgeries.

Our first result is the following.
We recall that all knots are classified into three families: torus knots, satellite knots, and hyperbolic knots.

\begin{theorem}\label{thm:main1}
There exist infinitely many pairs of non-equivalent knots $\{K_1,K_2\}$ in $S^3$
such that $m$-surgeries on them yield homeomorphic lens spaces for some integer $m$.
Moreover, $K_i$ can be chosen to be a torus knot, a satellite knot or
a hyperbolic knot, except that both of $K_1$ and $K_2$ cannot be satellite knots simultaneously. 
\end{theorem}

The exceptional case in Theorem \ref{thm:main1} is unavoidable as shown in Corollary \ref{cor:case-ss}, which
is obtained as a consequence of the next theorem.

\begin{theorem}\label{thm:main-torus}
\begin{enumerate}
\item There exist infinitely many pairs of non-equivalent torus knots in $S^3$
such that some half-integral surgeries on them yield homeomorphic lens spaces.
\item Let $K_1$ and $K_2$ be non-equivalent torus knots.
Suppose that a slope $r$ corresponds to a lens space surgery for both $K_1$ and $K_2$.
If the slope $r$ runs at least three times in the longitudinal direction, then $r$-surgeries on $K_1$ and $K_2$ cannot yield homeomorphic lens spaces.
\end{enumerate}
\end{theorem}

\begin{corollary}\label{cor:case-ss}
Non-equivalent satellite knots cannot yield homeomorphic lens spaces by the same Dehn surgery.
\end{corollary}

\begin{question}
Is there a lens space which can be obtained from three non-equivalent knots in $S^3$ by the same Dehn surgery?
\end{question}

Based on a computer experiment, we conjecture that the answer is negative.

The paper is organized as follows.
In Section \ref{sec:torus}, we give infinitely many pairs of torus knots that
yield homeomorphic lens spaces.
After establishing one result concerning a divisibility of integers by using
the classical theory of integral binary quadratic forms in Section \ref{sec:bqf},
we prove Theorem \ref{thm:main-torus} and Corollary \ref{cor:case-ss} in Section \ref{sec:non-int}. 
In Section \ref{sec:k+}, one special class of doubly primitive knots is reviewed.
In Section \ref{sec:hyperbolic}, infinitely many pairs of hyperbolic knots that yield
homeomorphic lens spaces are constructed by using tangles.
Finally, Section \ref{sec:diff} treats the case where the knots of a pair belong to different classes
to complete the proof of Theorem \ref{thm:main1}.

%%%%%%%%%%%%%%%%%%%%%%%%%%%%%%%%%%%%%%%%%%%%%%%%%%%%%%%%%%
\section{Torus knots}\label{sec:torus}

In this section, we give infinitely many pairs of torus knots which yield homeomorphic lens spaces by the same integral Dehn surgery.

Recall that the Fibonacci numbers are defined by the recurrence equation 
\[
F_{n+2}=F_{n+1}+F_{n}
\]
with $F_0=0$, $F_1=1$.
We make use of Cassini's identity (cf. \cite{GK})
\[
F_{k-1}F_{k+1}-F_k^2=(-1)^k
\]
for $k>0$.

Let $a_n=F_{n+2}$ and $b_n=F_{n+3}+F_{n+1}$ for $n\ge 1$.

\begin{lemma}\label{lem:torus1}
For any $n\ge 1$,
\[
a_{n+1}b_n+(-1)^{n+1}=a_nb_{n+1}+(-1)^n.
\]
\end{lemma}

\begin{proof}
By using Cassini's identity,
\begin{eqnarray*}
a_{n+1}b_n+(-1)^{n+1}&=& F_{n+3}(F_{n+3}+F_{n+1})+(-1)^{n+1}\\
                     &=& F_{n+3}^2+F_{n+3}F_{n+1}+(-1)^{n+1}\\
                     &=& F_{n+3}^2+F_{n+2}^2.
\end{eqnarray*}
Similarly,
\begin{eqnarray*}
a_{n}b_{n+1}+(-1)^{n}&=& F_{n+2}(F_{n+4}+F_{n+2})+(-1)^{n}\\
                     &=& F_{n+2}F_{n+4}+F_{n+2}^2+(-1)^{n}\\
                     &=& F_{n+3}^2+F_{n+2}^2.
\end{eqnarray*}
\end{proof}

As seen from Cassini's identity, two successive Fibonacci numbers are relatively prime.
Then it is easy to see that $\gcd(a_{n+1},b_n)=\gcd(a_{n},b_{n+1})=1$.

\begin{proposition}\label{prop:TT}
For $n\ge 1$, let $K$ be the torus knot of type $(a_{n+1},b_n)$, and $K'$ the torus knot of type $(a_n,b_{n+1})$.
Let $m=a_{n+1}b_n+(-1)^{n+1}\ (=a_nb_{n+1}+(-1)^n)$.
Then $K$ and $K'$ are not equivalent, and $m$-surgery on $K$ and $K'$ yield homeomorphic lens spaces.
\end{proposition}

\begin{proof}
Since $a_n<a_{n+1}<b_n<b_{n+1}$, $K$ and $K'$ are not equivalent.
By \cite{Mo}, $m$-surgery on $K$ and $K'$ yield the lens spaces $L(a_{n+1}b_n+(-1)^{n+1},a_{n+1}^2)$,
$L(a_{n}b_{n+1}+(-1)^{n},a_{n}^2)$, respectively.
Since $a_n^2+a_{n+1}^2=F_{n+2}^2+F_{n+3}^2=m$ as seen in the proof of Lemma \ref{lem:torus1},
$a_n^2+a_{n+1}^2\equiv 0 \pmod{m}$.
Thus these lens spaces are homeomorphic.
\end{proof}

%%%%%%%%%%%%%%%%%%%%%%%%%%%%%%%%%%%%%%%%%%%%%%%%
\section{Binary quadratic form}\label{sec:bqf}

In this section, we prove Proposition \ref{prop:bqf}, which will be used in Section \ref{sec:non-int}.
For its proof, we quickly review the classical theory of integral binary quadratic forms.
See \cite{F}, for example.

Let $f(x,y)=Ax^2+Bxy+Cy^2$ be an integral binary quadratic form with discriminant $\Delta=B^2-4AC$.
For our purpose, it is enough to assume that $\Delta$ is a positive nonsquare.
Let $m$ be a non-zero integer.
Then there is a finite algorithm to find all integral solutions $(x,y)\in \mathbb{Z}^2$
of the equation $f(x,y)=m$ as described below.

Let $\mathcal{S}=\{(x,y)\in \mathbb{Z}^2\mid f(x,y)=m\}$ be the set of integral solutions of $f(x,y)=m$.
Set 
\[
\rho=
\begin{cases}
\sqrt{\Delta/4} & \text{if $\Delta\equiv 0\pmod{4}$}, \\
\frac{1+\sqrt{\Delta}}{2} & \text{if $\Delta\equiv 1 \pmod{4}$}.\\
\end{cases}
\]
Let us consider the ring $\mathcal{O}_\Delta=\{x+y\rho \mid x,y\in \mathbb{Z}\}$.
Let $\mathcal{O}_{\Delta}^\times$ be the group of units of $\mathcal{O}_\Delta$, and
$\mathcal{O}_{\Delta,1}^\times=\{\alpha\in \mathcal{O}_\Delta^\times \mid N(\alpha)=1\}$
the subgroup of units for norm $1$.
Note that 
the norm $N(\alpha)$ of $\alpha=x+y\rho$ is given by the following:
\begin{equation*}
N(\alpha)=
\begin{cases}
x^2-\frac{\Delta}{4}y^2 & \text{if $\Delta\equiv 0\pmod{4}$},\\
x^2+xy-\frac{\Delta-1}{4}y^2 & \text{if $\Delta\equiv 1\pmod{4}$}.
\end{cases}
\end{equation*}
In fact, $\mathcal{O}_{\Delta,1}^\times$ corresponds to the solution set of the Pell equation $N(\alpha)=1$.
Then $\mathcal{O}_{\Delta,1}^\times$ acts on the set $\mathcal{S}$.
It is well known that the number of $\mathcal{O}_{\Delta,1}^\times$-orbits in $\mathcal{S}$ is finite.
Since $\mathcal{O}_{\Delta,1}^\times$ is infinite, the orbit of each solution is infinite, so
$\mathcal{S}$ is infinite, unless $\mathcal{S}=\emptyset$.
The action is explicitly given by the formulas
\begin{equation*}
(x',y')=
\begin{cases}
(x,y)\begin{pmatrix}
u-\frac{B}{2}v & Av \\
-Cv & u+\frac{B}{2}v \\
\end{pmatrix}
& \text{if $\Delta\equiv 0 \pmod{4}$}, \\
(x,y)\begin{pmatrix}
u+\frac{1-B}{2}v & Av \\
-Cv & u+\frac{1+B}{2}v \\
\end{pmatrix}
& \text{if $\Delta\equiv 1 \pmod{4}$},
\end{cases}
\end{equation*}
for $u+v\rho\in \mathcal{O}_{\Delta,1}^\times$ and $(x,y)\in \mathcal{S}$.

Let $\tau$ be the smallest unit of $\mathcal{O}_{\Delta,1}^\times$ that is greater than $1$.
Then every $\mathcal{O}_{\Delta,1}^\times$-orbit of integral solutions of $f(x,y)=m$ contains a solution $(x,y)\in \mathbb{Z}^2$
such that
\begin{equation*}
0\le y\le U=
\begin{cases}
|\frac{Am}{\Delta}(\tau+\bar{\tau}-2)|^{1/2} & \text{if $Am>0$},\\
|\frac{Am}{\Delta}(\tau+\bar{\tau}+2)|^{1/2} & \text{if $Am<0$},
\end{cases}
\end{equation*}
where $\bar{\tau}$ is the conjugate of $\tau$.
Furthermore, two distinct solutions $(x_1,y_1)$, $(x_2,y_2)\in \mathbb{Z}^2$ of $f(x,y)=m$
such that $0\le y_i\le U$ belong to the same $\mathcal{O}_{\Delta,1}^\times$-orbit if and only if
$y_1=y_2=0$ or $y_1=y_2=U$.

\begin{proposition}\label{prop:bqf}
Let $n\ge 3$ be an integer.
Let $a,b$ and $c$ be positive integers such that $a>1$ and $\gcd(a,b)=\gcd (a,c)=1$. 
Then $b^2\pm c^2$ is not divisible by $nabc\pm 1$.
\end{proposition}

\begin{proof}
Without loss of generality, we may assume that $b>c$.
Let $\varepsilon\in \{1,-1\}$.
If $b^2+c^2$ is divisible by $nabc+\varepsilon$,
then 
\begin{equation}\label{eq:start}
b^2+c^2=Q(nabc+\varepsilon)
\end{equation}
for some integer $Q\ge 1$.
Consider an integral binary quadratic form $f(x,y)=x^2-Qnaxy+y^2$.
Then the equation (\ref{eq:start}) means that
the equation $f(x,y)=\varepsilon Q$ has a solution $(b,c)$.

Similarly, if $b^2-c^2$ is divisible by $nabc+\varepsilon$,
then for a binary quadratic form $g(x,y)=x^2-Qnaxy-y^2$,
the equation $g(x,y)=\varepsilon Q$ has a solution $(b,c)$.
%%%
We remark that the discriminants $\Delta_f=(Qna)^2-4$ of $f$ and $\Delta_g=(Qna)^2+4$ of $g$ are positive and non-square.

First, we list all solutions in positive integers of the equation $f(x,y)=Q$.
For simplicity, let $\Delta=\Delta_f$.
Let $\mathcal{S}=\{(x,y)\in \mathbb{Z}^2\mid f(x,y)=Q\}$ be the set of all integral solutions of the
equation $f(x,y)=Q$.
Then the action of $\mathcal{O}_{\Delta,1}^\times$ on the set $\mathcal{S}$ is given by the formulas
\begin{equation}\label{eq:orbit}
(x',y')=
\begin{cases}
(x,y)\begin{pmatrix}
u+\frac{Qna}{2}v & v \\
-v & u-\frac{Qna}{2}v \\
\end{pmatrix}
& \text{if $\Delta\equiv 0 \pmod{4}$}, \\
(x,y)\begin{pmatrix}
u+\frac{1+Qna}{2}v & v \\
-v & u+\frac{1-Qna}{2}v \\
\end{pmatrix}
& \text{if $\Delta\equiv 1 \pmod{4}$},
\end{cases}
\end{equation}
for $u+v\rho\in \mathcal{O}_{\Delta,1}^\times$ and $(x,y)\in \mathcal{S}$.

Let $\tau$ be the smallest unit of $\mathcal{O}_{\Delta,1}^\times$ that is greater than $1$.
In fact, we see that 
\[
\tau=
\begin{cases}
\frac{Qna}{2}+\rho & \text{if $\Delta\equiv 0\pmod{4}$}, \\
\frac{Qna-1}{2}+\rho & \text{if $\Delta\equiv 1\pmod{4}$}.
\end{cases}
\]
Then every orbit contains a solution $(x,y)\in \mathbb{Z}^2$ such that
\[
0\le y\le U=\left|\frac{Q}{\Delta}(\tau+\bar\tau-2)\right|^{1/2}.
\]

In our case, $U<1$, and so $\mathcal{S}$ consists of a single $\mathcal{O}_{\Delta,1}^\times$-orbit.
Furthermore, $Q$ must be a square in order to $\mathcal{S}\ne \emptyset$. 
We start a solution $(\sqrt{Q},0)\in \mathcal{S}$.
By (\ref{eq:orbit}), 
\[
\tau\cdot (\sqrt{Q},0)=(\sqrt{Q},0)\begin{pmatrix}               
Qna & 1 \\
-1 & 0 
\end{pmatrix}
 =(Q^{3/2}na,\sqrt{Q}).
\]
Since
\[
(x,y)\begin{pmatrix}               
Qna & 1 \\
-1 & 0 
\end{pmatrix}
 =(Qnax-y,x),
\]
every solution in positive integers has a coordinate which is divisible by $a$.
Thus the equation $f(x,y)=Q$ cannot have the solution $(b,c)$, because $\gcd(a,b)=\gcd(a,c)=1$.
%%%%

For the equation $f(x,y)=-Q$, we have
\[
U=\left|\frac{-Q}{\Delta}(\tau+\bar\tau+2)\right|^{1/2}<1
\]
again.
However, $f(x,0)=x^2$ implies that the set of solutions of the equation $f(x,y)=-Q$ is empty.

Next, consider the equation $g(x,y)=Q$.
Let $\mathcal{T}=\{(x,y)\in \mathbb{Z}^2 \mid g(x,y)=Q\}$.
Put $\Delta=\Delta_g$.
Then $\mathcal{O}_\Delta$, $\mathcal{O}_{\Delta,1}^\times$ are defined in the same way, but
the action of $\mathcal{O}_{\Delta,1}^\times$ on the set $\mathcal{T}$ is given by the formulas
\begin{equation}\label{eq:orbit2}
(x',y')=
\begin{cases}
(x,y)\begin{pmatrix}
u+\frac{Qna}{2}v & v \\
v & u-\frac{Qna}{2}v \\
\end{pmatrix}
& \text{if $\Delta\equiv 0 \pmod{4}$}, \\
(x,y)\begin{pmatrix}
u+\frac{1+Qna}{2}v & v \\
v & u+\frac{1-Qna}{2}v \\
\end{pmatrix}
& \text{if $\Delta\equiv 1 \pmod{4}$},
\end{cases}
\end{equation}
for $u+v\rho\in \mathcal{O}_{\Delta,1}^\times$ and $(x,y)\in \mathcal{T}$.
Also, the smallest unit $\tau$ of $\mathcal{O}_{\Delta,1}^\times$ that is greater than $1$
is given by  
\[
\tau=
\begin{cases}
\left(\frac{Qna}{2}+\rho\right)^2=\frac{(Qna)^2}{2}+1+Qna\rho & \text{if $\Delta\equiv 0\pmod{4}$}, \\
\left(\frac{Qna-1}{2}+\rho\right)^2=\frac{(Qna)^2-Qna}{2}+1+Qna\rho & \text{if $\Delta\equiv 1\pmod{4}$}.
\end{cases}
\]

As before, we can evaluate
\[
U=\left|\frac{Q}{\Delta}(\tau+\bar\tau-2)\right|^{1/2}<\sqrt{Q}.
\]

On the other hand, if $(x,y)\in \mathcal{T}$, then
$\Delta y^2+4Q=(2x-Qnay)^2$.
That is, $\Delta y^2+4Q$ must be a square.
If $0<y<\sqrt{Q}$, then 
\[
Qnay<\sqrt{\Delta y^2+4Q}<Qnay+1.
\]
Hence $y=0$, and so $\mathcal{T}$ consists of a single $\mathcal{O}_{\Delta,1}^\times$-orbit.
Thus $Q$ must be a square in order to be $\mathcal{T}\ne \emptyset$.
Starting a solution $(\sqrt{Q},0)\in \mathcal{T}$,
\[
\tau\cdot (\sqrt{Q},0)=(\sqrt{Q},0)\begin{pmatrix}               
(Qna)^2+1 & Qna \\
Qna & 1 
\end{pmatrix}
 =(Q^{5/2}n^2a^2+\sqrt{Q},Q^{3/2}na)
\]
by the formulas (\ref{eq:orbit2}).
Thus for every solution in positive integers, the second coordinate is divisible by $a$.

Finally, for the equation $g(x,y)=-Q$, we have
\[
U=\left|\frac{-Q}{\Delta}(\tau+\bar\tau+2)\right|^{1/2}
\]
which is less than or equal to $\sqrt{Q}$ when $\Delta\equiv 0\pmod{4}$,
less than $\sqrt{Q}$ when $\Delta\equiv 1\pmod{4}$.

If $g(x,y)=-Q$, then $\Delta y^2-4Q=(2x-Qnay)^2$.
Thus $y\ne 0$. Furthermore,
if $y<\sqrt{Q}$, then
\[
Qnay-1<\sqrt{\Delta y^2-4Q}<Qnay.
\]
Therefore, $y=\sqrt{Q}$ is the only possibility, and so $Q$ must be a square.
As before, the set of solutions of the equation $g(x,y)=-Q$
consists of a single $\mathcal{O}_{\Delta,1}^\times$-orbit, whose representative is $(0,\sqrt{Q})$.
Then
\[
\tau\cdot (0,\sqrt{Q})=(0,\sqrt{Q})\begin{pmatrix}               
(Qna)^2+1 & Qna \\
Qna & 1 
\end{pmatrix}
 =(Q^{3/2}na,\sqrt{Q}).
\]
Hence the first coordinate is divisible by $a$ for any solution in positive integers.
\end{proof}

\begin{remark}
The requirement $a>1$ in the statement of Proposition \ref{prop:bqf} is necessary.
For example, let $a=1$, $b=3$, $c=8$ and $n=3$.
Then $b^2+c^2=73$ is divisible by $nabc+1=73$.
\end{remark}

%%%%%%%%%%%%%%%%%%%%%%%%%%%%%%%%%%%%%%%%%%%%%%%%%%%%%%%%%%%%%%%%%%%%%%%%%%%%%%%%%%%%%%%%%%%%%%%%%
\section{Non-integral surgery on torus knots}\label{sec:non-int}

In this section, we prove Theorem \ref{thm:main-torus}.

Let $\{a_n\}$ and $\{b_n\}$ be the sequences of positive integers defined by 

\begin{eqnarray}
a_{n+1} &=& a_n+b_{n}, \label{eqn:1} \\
b_{n+1} &=& a_{n+1}+a_n \label{eqn:2}
\end{eqnarray}
with $a_1=2$, $b_1=3$.

\begin{lemma}\label{lem:torus-half} For any $n\ge 1$,
\begin{enumerate}
\item $2a_nb_{n+1}+(-1)^{n+1}=2a_{n+1}b_n+(-1)^n$, 
\item $4a_{n+1}^2b_{n+1}^2+1=(2a_{n+1}b_{n+2}+(-1)^{n+2})(2a_{n}b_{n+1}+(-1)^{n+1})$.
\end{enumerate}
\end{lemma}

\begin{proof}
(1) By (\ref{eqn:1}) and (\ref{eqn:2}), $a_{n+1}=2a_{n}+a_{n-1}$.  Then
\begin{eqnarray*}
2a_nb_{n+1}-2a_{n+1}b_n &=& 2a_n(a_{n+1}+a_n)-2a_{n+1}(a_{n+1}-a_{n})\\
                        &=& 2(a_{n}^2-a_{n+1}^2+2a_na_{n+1})\\
                        &=& -2(a_{n-1}^2-a_n+2a_{n-1}a_n)\\
                        & & \vdots \\
                        &=& (-1)^{n-1}2(a_1^2-a_2^2+2a_1a_2)\\
                        &=& (-1)^n2 = (-1)^n-(-1)^{n+1}.
\end{eqnarray*}
(2) By (\ref{eqn:1}) and (\ref{eqn:2}), $2b_{n+1}=a_{n+2}+a_n$. 
Also, as shown above, $2a_nb_{n+1}-2a_{n+1}b_n=(-1)^n2$.
Thus, $a_nb_{n+1}-a_{n+1}b_n=(-1)^n$.
From (\ref{eqn:1}) and (\ref{eqn:2}),
$a_n(a_{n+1}+a_n)-a_{n+1}(a_{n+1}-a_n)=(-1)^n$.
Then $a_n^2+2a_na_{n+1}-a_{n+1}^2=(-1)^n$.
Thus,
\begin{eqnarray*}
\lefteqn{(2a_{n+1}b_{n+2}+(-1)^{n+2})(2a_{n}b_{n+1}+(-1)^{n+1})}\\
                        &=&(2a_{n+2}b_{n+1}+(-1)^{n+1})(2a_{n}b_{n+1}+(-1)^{n+1}) \\
                        &=& 4a_na_{n+2}b_{n+1}^2+(-1)^{n+1}2b_{n+1}(a_n+a_{n+2})+1 \\
                        &=& 4b_{n+1}^2(a_na_{n+2}+(-1)^{n+1})+1 \\
                        &=& 4b_{n+1}^2(a_n(a_{n+1}+b_{n+1})+(-1)^{n+1})+1 \\
                        &=& 4b_{n+1}^2(a_na_{n+1}+a_n(a_{n+1}+a_n)+(-1)^{n+1})+1 \\
                        &=& 4b_{n+1}^2(2a_na_{n+1}+a_n^2+(-1)^{n+1})+1 \\
                        &=& 4b_{n+1}^2a_{n+1}^2+1.
\end{eqnarray*}
\end{proof}

From Lemma \ref{lem:torus-half}(1),
we have that $\gcd(a_n,b_{n+1})=\gcd(b_n,a_{n+1})=1$.

\begin{proof}[Proof of Theorem \textup{\ref{thm:main-torus}(1)}]
Let $K_1$ be the torus knot of type $(a_n,b_{n+1})$, $K_2$ the torus knot of type $(b_n,a_{n+1})$.
Since $a_n<b_n<a_{n+1}<b_{n+1}$ for any $n\ge 1$,
$K_1$ and $K_2$ are not equivalent.
Then $\frac{2a_nb_{n+1}+(-1)^{n+1}}{2}$-surgery on $K_1$ and
$\frac{2a_{n+1}b_n+(-1)^n}{2}$-surgery on $K_2$ yield the lens spaces $L(2a_nb_{n+1}+(-1)^{n+1},2b_{n+1}^2)$ and
$L(2a_{n+1}b_n+(-1)^n,2a_{n+1}^2)$, respectively.
By Lemma \ref{lem:torus-half}, the surgery coefficients are the same, and the two lens spaces are homeomorphic.
\end{proof}

%%%%
In the rest of this section, we prove (2) of Theorem \ref{thm:main-torus}, and give a proof of Corollary \ref{cor:case-ss}.

Let $K_1$ be the torus knot of type $(p,q)$, $K_2$ the torus knot of type $(r,s)$.
Suppose $n\ge 3$.
If $m/n$-surgery on $K_1$ yields a lens space, then $\Delta(pq/1,m/n)=|npq-m|=1$, so $m=npq\pm 1$.
Hence if $m/n$-surgery on $K_1$ and $K_2$ yield homeomorphic lens spaces, then
$npq+\varepsilon=nrs+\varepsilon'$ for some $\varepsilon, \varepsilon'\in \{1,-1\}$.
Since we consider non-trivial torus knots, we can assume that $p,q,r$ and $s$ are positive
by taking mirror images, if necessary.
Moreover, we may assume that $2\le q<p$, $2\le s<r$ and $r<p$.
Because $n\ge 3$, $\varepsilon=\varepsilon'$, and so $pq=rs$.
By \cite{Mo}, $m/n$-surgery on $K_1$ and $K_2$ yield $L(m,nq^2)$ and $L(m,ns^2)$, respectively.

Theorem \ref{thm:main-torus}(2) follows directly from the following.

\begin{proposition}
Two lens spaces $L(m,nq^2)$ and $L(m,ns^2)$ are not homeomorphic.
\end{proposition}

\begin{proof}
The two lens spaces are homeomorphic if and only if 
\begin{eqnarray}
nq^2 &\equiv& \pm ns^2 \pmod{m}\ \text{or} \label{eq:ss1}\\
n^2q^2s^2 &\equiv& \pm 1 \pmod{m} \label{eq:ss2}.
\end{eqnarray}
Since $npq+\varepsilon=nrs+\varepsilon$, $npq=nrs$.
Thus $q<s<r<p$.

First, $nq^2\not\equiv ns^2\pmod{m}$, because
$0<n(s^2-q^2)<ns^2<nrs-1\le m$.
Next, since $nq^2+ns^2<n(pq+rs)-2=2npq-2\le 2m$,
the equation $nq^2\equiv -ns^2\pmod{m}$ is possible only when $nq^2+ns^2=m$.
However, this is impossible, because $m$ is not divisible by $n$.

The impossibility of the equation (\ref{eq:ss2}) will be shown in the next proposition.
\end{proof}

\begin{proposition}
$n^2q^2s^2 \not\equiv \pm 1 \pmod{m}$.
\end{proposition}

\begin{proof}
Suppose $n^2q^2s^2\equiv 1 \pmod{m}$.
Then $n^2q^2s^2-1=km$ for some integer $k\ge 1$.
Since $-1\equiv k\varepsilon \pmod{n}$ (recall $m=npq+\varepsilon$),
$k\equiv -\varepsilon \pmod{n}$.
Put $k=n\ell-\varepsilon$ with $\ell\ge 1$.
(If $\ell=0$, then $k=-\varepsilon=-1$, so $n^2q^2s^2-1=m$.
This implies that $q$ divides $p$, a contradiction.)
Then $n^2q^2s^2-1=(n\ell-\varepsilon)(npq+\varepsilon)$ implies
\begin{equation}\label{eq:ell}
q(nqs^2-p(n\ell-\varepsilon))=\varepsilon\ell.
\end{equation}
Thus $q$ divides $\ell$, and $\gcd(p,s)$ divides $\ell/q$.
For simplicity, we denote $\gcd(x,y)$ by $(x,y)$.

Hence
\begin{equation}\label{eq:gcd}
\begin{split}
p(n\ell-\varepsilon) &= nqs^2-\frac{\varepsilon\ell}{q}\\
                     &= nq (p,s)^2(q,s)^2-\frac{\varepsilon \ell}{q}\\
                     &= (p,s)\left(\frac{nq}{(q,s)}(q,s)^3(p,s)-\frac{\varepsilon\ell}{q (p,s)}\right)
\end{split}
\end{equation}
Here we put $a=\frac{q(p,s)}{(q,s)}$, $b=(q,s)$ and $c=\frac{\ell}{q(p,s)}$.
Then $abc=\ell$.

\begin{claim}\label{cl:a1}
$a>1$.
\end{claim}

\begin{proof}[Proof of Claim \ref{cl:a1}]
Assume $a=1$.
Then $(p,s)=1$ and $q=(q,s)$. Since $s=(p,s)(q,s)$, $s=(q,s)$.
Thus $s=q$, so $p=r$, a contradiction.
\end{proof}

\begin{claim}\label{cl:abc}
$(a,b)=(a,c)=1$.
\end{claim}

\begin{proof}[Proof of Claim \ref{cl:abc}]
First, $(p,s)$ and $(q,s)$ are coprime.
Also, $q/(q,s)$ and $(q,s)$ are coprime, otherwise $(r,s)>1$.
Thus $(a,b)=1$.

Next, assume $(a,c)>1$.
Let $d$ be a prime factor of $(a,c)$.
From the equation (\ref{eq:ell}),
\[
nqs^2-p(n\ell-\varepsilon)=\varepsilon \frac{\ell}{q}.
\]
Dividing it by $(p,s)$ gives
\begin{equation}\label{eq:dd}
nqs\frac{s}{(p,s)}-\frac{p}{(p,s)}(n\ell-\varepsilon)=\varepsilon c.
\end{equation}
Since $d$ divides $a$, $d$ divides $q$ or $s$.
Similarly, $d$ divides $\ell$, since $d$ divides $c$. 
Thus the equation (\ref{eq:dd}) gives
\[
\frac{p}{(p,s)}\varepsilon\equiv 0 \pmod{d}.
\]
However, this is impossible, because
$(p,s)$ and $p/(p,s)$ are coprime.
\end{proof}

On the other hand, the equation (\ref{eq:gcd}) yields
\[
p(nabc-\varepsilon)=(p,s)(nab^3-\varepsilon c).
\]
Because $(p,s)$ divides $p$, this equation means that
$nab^3-\varepsilon c$ is divisible by $nabc-\varepsilon$.
Furthermore,
\[
\frac{nab^3-\varepsilon c}{nabc-\varepsilon}=c+\frac{nab(b^2-c^2)}{nabc-\varepsilon}
\]
implies that $b^2-c^2$ is divisible by $nabc-\varepsilon$, since $nab$ and $nabc-\varepsilon$ are coprime.

Similarly, if $n^2q^2s^2\equiv -1 \pmod{m}$, then we have the fact that $b^2+c^2$ is divisible by $nabc-\varepsilon$.
However, these are impossible by Proposition \ref{prop:bqf}.
\end{proof}

\begin{proof}[Proof of Corollary \ref{cor:case-ss}]
Among satellite knots,
only the $(2,2pq+\varepsilon)$-cable $K$ of the $(p,q)$-torus knot admits a lens space surgery for $\varepsilon=\pm 1$.
Then the slope is $4pq+\varepsilon$, and $L(4pq+\varepsilon,4q^2)$ arises.
This surgery on $K$ is equivalent to $\frac{4pq+\varepsilon}{4}$-surgery on its companion torus knot.
Thus the result follows from Theorem \ref{thm:main-torus}(2).
\end{proof}

%%%%%%%%%%%%%%%%%%%%%%%%%%%%%%%%%%%%%%%%%%%%%%%%%%%%%%%%%%%%%%%%%%%%%%%%%%
\section{Doubly primitive knot}\label{sec:k+}

In this section, we study a special class of doubly primitive knots $k^+(a,b)$ defined by Berge \cite{B}.
In particular, two infinite sequences of $k^+(a,b)$ are proved to be hyperbolic via dual knots in lens spaces.
As far as we know,  the determination of hyperbolicity of $k^+(a,b)$ is still open, in general.

\begin{figure}
\centering
\includegraphics[scale=0.45]{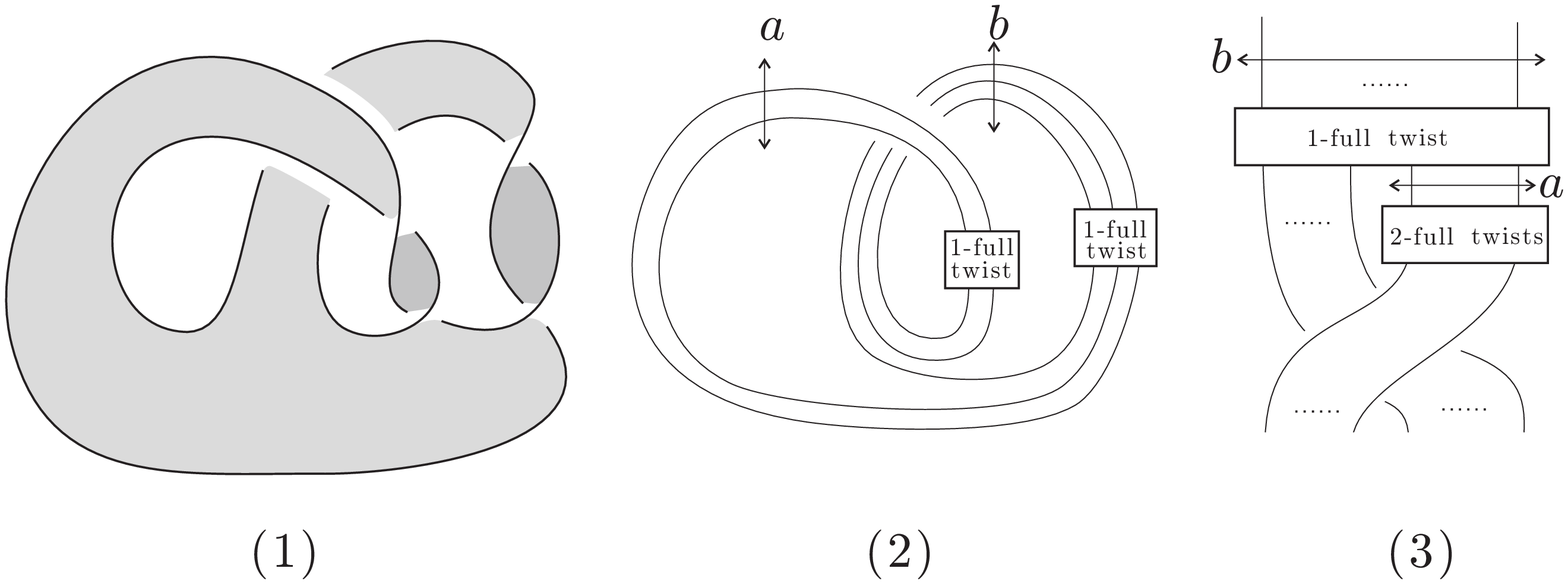}
\caption{$k^+(a,b)$}
\label{fig:trefoil}
\end{figure}

For a pair $(a,b)$ of coprime positive integers,
let $k^+(a,b)$ denote the doubly primitive knot defined by Berge \cite{B}, which
lies on a genus one fiber surface of the left-handed trefoil as shown in Figure \ref{fig:trefoil}(1).
Then $(a^2+ab+b^2)$-surgery on $k^+(a,b)$ yields the lens space $L(a^2+ab+b^2,(a/b)^2)$,
where $a/b$ is calculated in $\mathbb{Z}_{a^2+ab+b^2}$.
(We adopt the notation in \cite{Y}, but the orientation of lens spaces is opposite to ours). 
We remark that $k^+(a,b)$ and $k^+(b,a)$ are equivalent by the symmetry of the fiber surface.
For example, $k^+(1,3)$ is the $(3,4)$-torus knot whose $13$-surgery yields $L(13,9)$, and 
$k^+(2,3)$, as shown in Figure \ref{fig:trefoil}(2), is the $(-2,3,7)$-pretzel knot whose $19$-surgery
yields $L(19,7)$.

\begin{lemma}\label{lem:k+genus}
$k^+(a,b)$ is a fibered knot with genus $\frac{(a+b-1)^2-ab}{2}$.
\end{lemma}

\begin{proof}
It is easy to see that $k^+(a,b)$ has a form of the closure of a positive braid as shown in Figure \ref{fig:trefoil}(3).
By \cite{St}, Seifert's algorithm gives a fiber surface.
The braid has $b$ strings and $a^2+ab+b^2-2a-b$ crossings, so
the fiber has the given genus.
(See also \cite[Corollary 3]{Y} or \cite{HM}.)
\end{proof}

%%%%%%%%%%%%%%%%%%%%%%%%%%%%%%%%%%%%%%%%%%%%
In general, let $K$ be a knot in $S^3$ whose $p$-surgery yields $L(p,q)$ with $p>q>0$.
Then the core $K^*$ of the attached solid torus of $K(p)$ is called the \textit{dual knot\/} of $K$ (with respect to $p$-surgery).
Berge \cite{B} shows that 
if $K$ is a doubly primitive knot whose surface slope is $p$, then
$K^*$ is a $(1,1)$-knot in $L(p,q)$, and it has a canonical form parameterized by 
a single integer $k$ with $0< k<p$ (see \cite{Sa0,Sa1}).
Following \cite{Sa0}, we denote it by $K(L(p,q);k)$.
It is known that $K(L(p,q);k)$ is isotopic to $K(L(p,q);p-k)$.

For $n=1,2,\dots,p-1$,
let $\phi_n$ be an integer such that $\phi_n\equiv nq \pmod{p}$ and $0<\phi_n<p$.
We call this finite sequence $\{\phi_n\}$ the \textit{basic sequence\/} for $(p,q)$.
Because of $\gcd(p,q)=1$, $\phi_n$'s are mutually distinct.
In particular, $k$ appears in the basic sequence.
Let $h$ be the position of $k$, that is, $\phi_h=k$.
Here, set
\begin{eqnarray*}
s&=&\sharp \{i\mid i<k\ \text{and $i$ appears before $k$ in the basic sequence} \},\\
\ell&=& \sharp \{i\mid i>k\ \text{and $i$ appears before $k$ in the basic sequence} \},\\
s'&=&\sharp \{i\mid i<k\ \text{and $i$ appears after $k$ in the basic sequence} \}, \\
\ell'&=& \sharp \{i\mid i>k\ \text{and $i$ appears after $k$ in the basic sequence} \}.\\
\end{eqnarray*}
%Then $s+\ell=h-1$, $s+s'=k-1$ and $\ell+\ell'=p-1-k$.
Let
\[
\Phi=\min\{s,s',\ell,\ell'\}.
\]
This is determined for the triplet $(p,q,k)$, so for the dual knot $K(L(p,q);k)$.
However, the main result of \cite{Sa} says that $\Phi$ depends only on the original knot $K$ and
a lens space surgery slope $p$, and that 
$K$ is hyperbolic if and only if $\Phi\ge 2$, equivalently,
each of $s,s',\ell,\ell'$ is at least two.

For $k^+(a,b)$, let $p=a^2+ab+b^2$.
Then $p$-surgery yields a lens space $L(p,q)$ where $q\equiv (\frac{b}{a+b})^2$.
(Note that $(\frac{a}{b})^2\equiv (\frac{b}{a+b})^2 \pmod{p}$.)
By \cite{Sa0}, the dual knot is represented as $K(L(p,q);k)$ with
$k\equiv -\frac{b}{a+b}\pmod{p}$. 
(By definition, the parameter $k$ is chosen so as $0<k<p$.)

\begin{lemma}\label{lem:kq-relation}
Let $p$, $q$ and $k$ be defined as above.
\begin{enumerate}
\item $k+q+1\equiv 0 \pmod{p}$.
\item $k\equiv q^2 \pmod{p}$.
\item $kq\equiv 1 \pmod{p}$.
\end{enumerate}
\end{lemma}

\begin{proof}
(1) 
\begin{eqnarray*}
k+q+1&\equiv&-\frac{b}{a+b}+\frac{b^2}{(a+b)^2}+1=\frac{-b(a+b)+b^2+(a+b)^2}{(a+b)^2}\\
     &=& \frac{(a+b)^2-ab}{(a+b)^2}=\frac{p}{(a+b)^2}\equiv  0 \pmod{p}.
\end{eqnarray*}

(2) By (1), $q^2-k\equiv (-k-1)^2-k=k^2+k+1\equiv q+k+1\equiv 0 \pmod{p}$.

(3) Similarly, $kq\equiv k(-k-1)=-k^2-k\equiv -q-k\equiv 1 \pmod{p}$ by (1).
\end{proof}

%%%%%%%%%
\subsection{$k^+(3n+1,3n+4)$}

For $k^+(3n+1,3n+4)$, let $p=27n^2+45n+21$.
Then $p$-surgery yields a lens space $L(p,q)$ with $q=(3n+2)^2$, and 
the dual knot is $K(L(p,q);k)$ with $k\equiv -(3n+2)^2-1 \pmod{p}$.

\begin{lemma}\label{lem:k+1}
For $n\ge 1$, $k^+(3n+1,3n+4)$ is hyperbolic.
\end{lemma}

\begin{proof}
Let $a=3n+1$, $b=3n+4$ and $k_0=p-q-1$. Then 
direct calculations show that $3q<p<4q$,
$k_0\equiv k \pmod{p}$ and $2q-1<k_0<3q-a$.
Thus we can use the triplet $(p,q,k_0)$ to calculate the invariant $\Phi$.

Let $\{\phi_i\}$ be the basic sequence.
(Recall that any term $\phi_i$ of the basic sequence is chosen so that $0<\phi_i<p$.)
Since $q^2\equiv k\equiv k_0\pmod{p}$ by Lemma \ref{lem:kq-relation}, $\phi_q=k_0$. 

First, we study the four consecutive terms $\phi_{a+b},\phi_{a+b+1},\phi_{a+b+2},\phi_{a+b+3}$,
which appear before $k_0$.
Since $(a+b)q\equiv p-a \pmod{p}$, $\phi_{a+b}=p-a$.
Then
\[
q-a<2q-a<k_0<3q-a<p-a,
\]
so $\phi_{a+b+1}=q-a$, $\phi_{a+b+2}=2q-a$, $\phi_{a+b+3}=3q-a$.
Hence
$\phi_{a+b}>k_0$, $\phi_{a+b+1}<k_0$, $\phi_{a+b+2}<k_0$ and $\phi_{a+b+3}>k_0$.

Similarly, we study the four consecutive terms right after $k_0$.
(Since $q+4<p-1$, there are more than four terms after $k_0$.)
Since $k_0+q=p-1$, $\phi_{q+1}=p-1$.
Then $\phi_{q+2}=q-1$,  $\phi_{q+3}=2q-1$,  $\phi_{q+4}=3q-1$.
Hence $\phi_{q+1}>k_0$, $\phi_{q+2}<k_0$, $\phi_{q+3}<k_0$ and $\phi_{q+4}>k_0$.
Thus we have $\Phi\ge 2$, showing that the dual knot (and the original knot) is hyperbolic.
\end{proof}

%%%%%%%%%%%%%%%%%%%%%%%%
\subsection{$k^+(F_{n+2},F_n)$}

For Fibonacci numbers, see Section \ref{sec:torus}.
Let $p=F_{n}^2+F_{n}F_{n+2}+F_{n+2}^2$.
Then $p$-surgery on $k^+(F_{n+2},F_n)$ yields a lens space $L(p,q)$ with $q\equiv (\frac{F_n}{F_{n+2}+F_n})^2 \pmod{p}$,
and let the dual knot denote by $K(L(p,q);k)$.

\begin{lemma}\label{lem:F2}
$p=4F_nF_{n+2}+(-1)^n$, $q\equiv (-1)^{n+1}4F_n^2 \pmod{p}$, and
$k\equiv (-1)^n4F_n(F_n+F_{n+2}) \pmod{p}$.
\end{lemma}

\begin{proof}
By Cassini's identity $F_nF_{n+2}-F_{n+1}^2=(-1)^{n+1}$,
\begin{eqnarray*}
4F_nF_{n+2}+(-1)^n &=& 4F_{n+1}^2+3(-1)^{n+1}\\
                   &=& 3(F_{n+1}^2+(-1)^{n+1})+F_{n+1}^2\\
                   &=& 3F_nF_{n+2}+F_{n+1}^2\\
                   &=& F_nF_{n+2}+2F_n(F_{n}+F_{n+1})+F_{n+1}^2\\
                   &=& F_nF_{n+2}+2F_n^2+2F_nF_{n+1}+F_{n+1}^2\\
                   &=& F_nF_{n+2}+F_n^2+(F_n+F_{n+1})^2\\
                   &=& F_nF_{n+2}+F_n^2+F_{n+2}^2=p.
\end{eqnarray*}

Thus $4F_nF_{n+2}+(-1)^n\equiv 0 \pmod{p}$.
To show $q\equiv (-1)^{n+1}4F_n^2 \pmod{p}$,
it suffices to show $(-1)^{n+1}4(F_n+F_{n+2})^2\equiv 1\pmod{p}$.
This follows from the equation $(F_n+F_{n+2})^2\equiv F_nF_{n+2} \pmod{p}$.

Finally, 
$(-1)^n4F_n(F_n+F_{n+2})=(-1)^n4(F_n^2+F_nF_{n+2})\equiv (-1)^{n+1}4F_{n+2}^2 \pmod{p}$.
Then $(-1)^{n+1}4F_{n+2}^2 q\equiv (4F_{n+2}F_{n})^2\equiv 1\pmod{p}$.
This shows that $(-1)^{n+1}4F_{n}(F_n+F_{n+2})\equiv 1/q\equiv k \pmod{p}$ by Lemma \ref{lem:kq-relation}(3).
\end{proof}

\begin{lemma}\label{lem:k+2}
For $n\ge 3$, $k^+(F_{n+2},F_{n})$ is hyperbolic.
\end{lemma}

\begin{proof}
As mentioned above, $p$-surgery on $k^+(F_{n+2},F_{n})$ yields $L(p,q)$.
Consider the dual knot $K(L(p,q);k)$ in $L(p,q)$.

%%%%%%%%%%%%%%%%%%%%%%%%%%%%%%%%%%%%%%%%
First, we assume that $n$ is odd.
Then $p=4F_nF_{n+2}-1$, $q\equiv 4F_n^2 \pmod{p}$, and
$k\equiv -4F_n(F_n+F_{n+2}) \pmod{p}$ by Lemma \ref{lem:F2}.

To make the calculation of the invariant $\Phi$ easier,
put $q_0=p-4F_n^2$ and $k_0=p-q_0+1$.
Then $0<q_0<p$, $0<k_0<p$, and $q_0\equiv -q\pmod{p}$ and $k_0\equiv -k \pmod{p}$.

\begin{claim}\label{cl:F2-1}
$3q_0/2<p<2q_0$ and $2q_0-p<k_0<q_0$.
\end{claim}

\begin{proof}[Proof of Claim \ref{cl:F2-1}]
$2q_0-p=p-8F_n^2=4F_nF_{n+2}-1-8F_n^2=4F_n(F_n+F_{n+1})-8F_n^2-1=4F_n(F_{n+1}-F_n)-1\ge 7$.
Since $3F_n>F_{n+2}$, $2p-3q_0=12F_n^2-p=4F_n(3F_n-F_{n+2})+1\ge 9$.

Next, $k_0-2q_0+p=2p-3q_0+1\ge 10$.
Finally, $q_0-k_0=2q_0-p-1\ge 6$.
\end{proof}

%%%%
For $(p,q_0)$, let $\{\phi_i\}$ be the basic sequence.
Let $h=p-q_0$.
Since $hq_0\equiv k_0 \pmod{p}$, the number $k_0$ appears in the sequence as the $h$-th term.
Note that $h>4$, because $2p-3q_0\ge 9$.

To evaluate $\Phi$, we investigate some specific terms in the basic sequence. 
We have $\phi_1=q_0>k_0$ and $\phi_2=2q_0-p<k_0$.
Also, $\phi_{h-1}=k_0-q_0+p>k_0$ and $\phi_{h-2}=k_0-2q_0+p<k_0$.
Since $h>4$, these four terms $\phi_1,\phi_2,\phi_{h-2},\phi_{h-1}$ are distinct.
Next, $\phi_{p-1}=p-q_0<k_0$ and $\phi_{p-2}=2p-2q_0>k_0$.
Since $2h<p$, $h<p-h+1$.
Thus $\phi_{p-h+1}$ and $\phi_{p-h+2}$, which are distinct from $\phi_{p-1}$ and $\phi_{p-2}$, appear after $k_0$ in the basic sequence.
We have
$\phi_{p-h+1}=q_0-k_0$, since $(p-h+1)q_0\equiv q_0-k_0\pmod{p}$.
Then $k_0-(q_0-k_0)=2k_0-q_0=2p-3q_0+2>0$ implies $\phi_{p-h+1}<k_0$.
Finally, $\phi_{p-h+2}=2q_0-k_0>k_0$.
Again, the fact $h>4$ means that the four terms $\phi_{p-h+1},\phi_{p-h+2},\phi_{p-2},\phi_{p-1}$ are distinct.
Hence $\Phi\ge 2$.

%%%%%%%%%%%%%%%%%%%%%%%%%%%%%%%%%%%%%%%
Second, assume that $n$ is even.
Then $p=4F_nF_{n+2}+1$, $q\equiv -4F_n^2 \pmod{p}$, and
$k\equiv 4F_n(F_n+F_{n+2}) \pmod{p}$ by Lemma \ref{lem:F2}.
In this case, put $q_0=p-4F_n^2+1$ and $k_0=p-q_0+1$.
Then $0<q_0<p$ and $0<k_0<p$.
It is easy to check that Claim \ref{cl:F2-1} holds without any change.

By Lemma \ref{lem:kq-relation}, $q_0q\equiv (q+1)q\equiv k+q\equiv -1 \pmod{p}$ and $k_0\equiv -q\pmod{p}$.
Under a (orientation-reversing) homeomorphism from $L(p,q)$ to $L(p,q_0)$,
the dual knot $K(L(p,q);k)$ is mapped to $K(L(p,q_0);k_0)$ (see \cite{Sa}).
Thus we can use $(p,q_0,k_0)$, instead of $(p,q,k)$, to evaluate $\Phi$.

By Lemma \ref{lem:kq-relation}(2), $q^2\equiv k \pmod{p}$.
Thus $16F_n^4\equiv 4F_n(F_n+F_{n+2})$.
Hence $q_0^2+k_0 \equiv (1-4F_n^2)^2+4F_n^2\equiv 16F_n^4-4F_n^2+1
\equiv 4F_n(F_n+F_{n+2})-4F_n^2+1\equiv 4F_n^2F_{n+2}^2+1 \equiv 0 \pmod{p}$.
This means that $(p-q_0)q_0\equiv k_0 \pmod{p}$.
Let $h=p-q_0$.
Then, $k_0$ appears in the basic sequence for $(p,q_0)$ as the $h$-term.
Since $h>4$, the argument in the case where $n$ is odd works verbatim, so we have $\Phi\ge 2$.
\end{proof}

%%%%%%%%%%%%%%%%%%%%%%%%%%%%%%%%%%%%%%%%%
\section{Hyperbolic knots}\label{sec:hyperbolic}

A Seifert fibered manifold is said to be \textit{of type\/} $X(p_1,p_2,\dots,p_n)$ if
it admits a Seifert fibration over the surface $X$ with $n$ exceptional fibers of indices $p_1,p_2,\dots,p_n$.
In this paper, $X$ will be either the $2$-sphere $S^2$ or the disk $D^2$.

\begin{figure}
\centering
\includegraphics[scale=0.5]{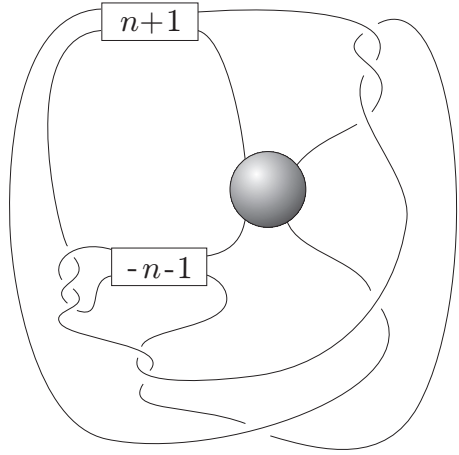}
\caption{The tangle $B_n$}
\label{fig:caseHH}
\end{figure}

For $n\ge 1$, let $B_n$ be the tangle as shown in Figure \ref{fig:caseHH},
where a rectangle denotes horizontal half-twists.
If the number is positive, the twist is right-handed, otherwise, left-handed.

\begin{figure}
\centering
\includegraphics[scale=0.45]{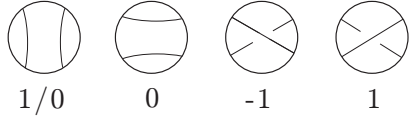}
\caption{Some rational tangles}
\label{fig:ratio}
\end{figure}

Given $\alpha\in \mathbb{Q}\cup \{1/0\}$,
$B_n(\alpha)$ denotes the knot or link in $S^3$ obtained by inserting the rational tangle of slope $\alpha$
into the central puncture of $B_n$.
Also, $\widetilde{B}_n$ is the double branched cover of $S^3$ branched over $B_n(\alpha)$. 
In fact, we need only four rational tangles as shown in Figure \ref{fig:ratio}.

\begin{lemma}\label{lem:caseHH-tangle}
\begin{enumerate}
\item $\widetilde{B}_n(1/0)=S^3$.
\item $\widetilde{B}_n(0)=L(27n^2+45n+21,-9n^2-12n-5)$.
\item $\widetilde{B}_n(1)$ is a Seifert fibered manifold of type $S^2(2,n+2,15n+11)$.
\item $\widetilde{B}_n(-1)$ is a non-Seifert toroidal manifold $D^2(2,n)\cup D^2(2,3n+1)$, which
contains a unique incompressible torus, if $n\ge 2$, or
a Seifert fibered manifold of type $S^2(2,3,4)$ if $n=1$.
\end{enumerate}
\end{lemma}

\begin{figure}
\centering
\includegraphics[scale=0.4]{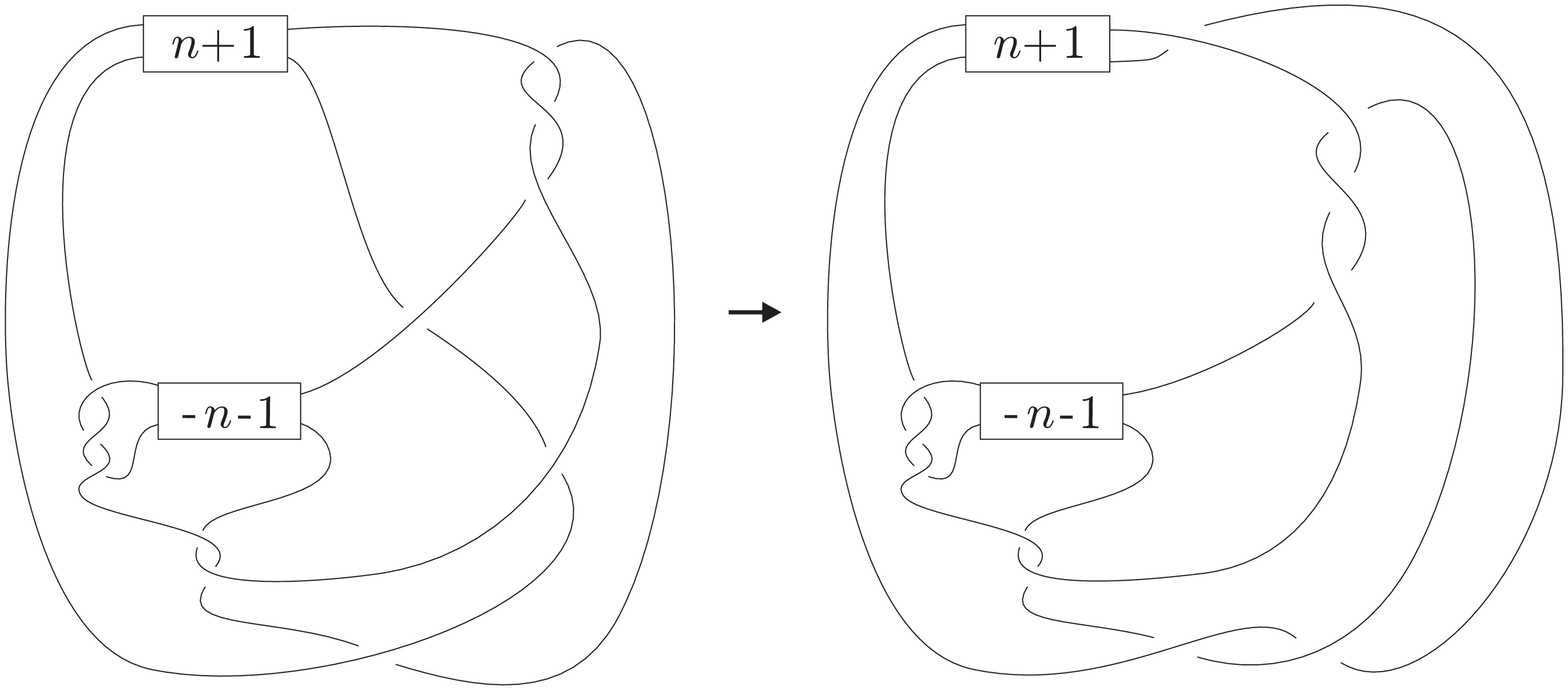}
\caption{$B_n(1)$}
\label{fig:HH-1}
\end{figure}

\begin{figure}
\centering
\includegraphics[scale=0.4]{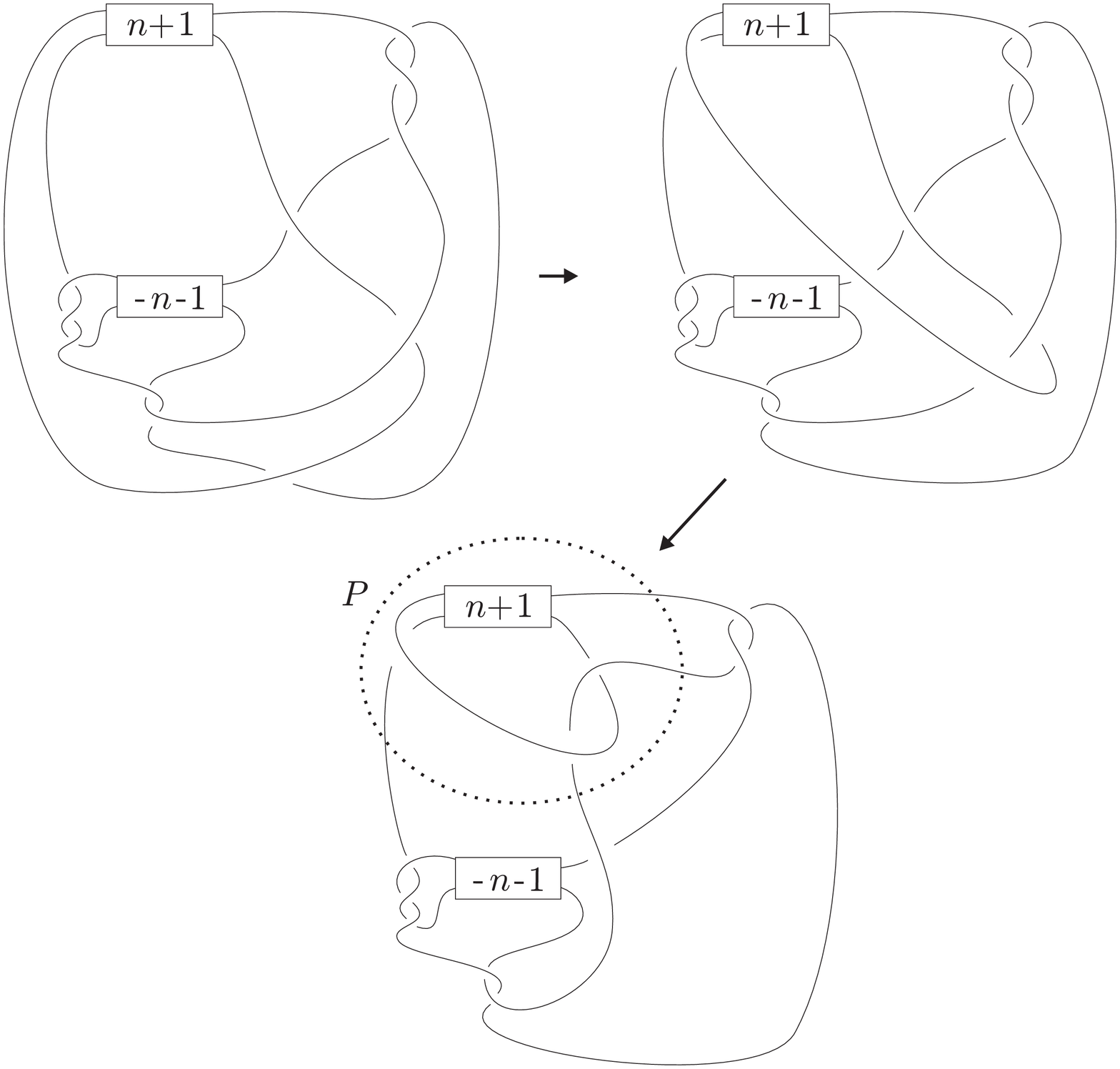}
\caption{$B_n(-1)$}
\label{fig:HH1}
\end{figure}

\begin{proof}
It is straightforward to check that $B_n(1/0)$ is the unknot and that
$B_n(0)$ is the $2$-bridge knot corresponding to $-(9n^2+12n+5)/(27n^2+45n+21)$.

Figure \ref{fig:HH-1} shows that $B_n(1)$ is a Montesinos link or knot of length three.
Thus $\widetilde{B}_n(1)$ is a Seifert fibered manifold of type $S^2(2,n+2,15n+11)$.

Figure \ref{fig:HH1} shows that $B_n(-1)$ is decomposed along a tangle sphere $P$ into
two tangles.
If $n>1$, then each side of $P$ is a Montesinos tangle.
Thus $\widetilde{B}_n(1)$ is decomposed along a torus into two Seifert fibered manifolds over the disk with
two exceptional fibers.
Since Seifert fibers on both sides intersect once on the torus, 
$\widetilde{B}_n(-1)$ is not Seifert.
It is well known that such a $3$-manifold contains a unique incompressible torus.
When $n=1$, $B_n(-1)$ is a Montesinos link of length three.
Hence $\widetilde{B}_n(-1)$ is a Seifert fibered manifold over the $2$-sphere with three exceptional fibers.
\end{proof}

By Lemma \ref{lem:caseHH-tangle}(1), the lift of $B_n$ in $\widetilde{B}_n(1/0)$ gives
the knot exterior of some knot $K_n$ in $S^3$, which is uniquely determined by Gordon-Luecke's theorem \cite{GL}.
Furthermore, $K_n$ admits integral Dehn surgeries yielding a lens space, a Seifert fibered manifold, and a toroidal manifold (unless $n=1$) by
Lemma \ref{lem:caseHH-tangle}.

The following criterion of hyperbolicity is used also in Section \ref{sec:diff}.

\begin{lemma}\label{lem:criterionH}
Let $K$ be a knot in $S^3$.
If $K$ admits an integral lens space surgery $m$, and
neither $K(m-1)$ nor $K(m+1)$ has a lens space summand, then $K$ is hyperbolic.
\end{lemma}

\begin{proof}
Assume not.  Then $K$ is either a torus knot or a satellite knot.
For the (non-trivial) $(p,q)$-torus knot, the only integral lens space surgery slopes are $pq-1$ and $pq+1$, and
$pq$-surgery yields the connected sum of two lens spaces by \cite{Mo}.
Thus $K$ is not a torus knot.

Assume that $K$ is a satellite knot.
Since $K$ has a lens space surgery, $K$ is the $(2,2pq+\varepsilon)$-cable of the $(p,q)$-torus knot where $\varepsilon\in \{1,-1\}$
by \cite{BL,W,Wu}.
Then the lens space surgery is $4pq+\varepsilon$.
However, the adjacent slope $4pq+2\varepsilon$ is equal to the cabling slope, and so
$K(4pq+2\varepsilon)$ has a lens space summand, a contradiction.
Thus $K$ is hyperbolic.
\end{proof}

\begin{lemma}\label{lem:hyp1}
$K_n$ is hyperbolic.
\end{lemma}

\begin{proof}
This immediately follows from
Lemmas \ref{lem:caseHH-tangle} and \ref{lem:criterionH}.
\end{proof}

\begin{lemma}\label{lem:caseHH-Hlens}
The knot $K_n$ defined above satisfies the following.
\begin{enumerate}
\item The genus of $K_n$ is $(27n^2+33n+10)/2$.
\item Let $m=27n^2+45n+21$.
Then $m$-surgery on $K$ yields the lens space $L(m,-9n^2-12n-5)$. 
\end{enumerate}
\end{lemma}

\begin{figure}
\centering
\includegraphics[scale=0.4]{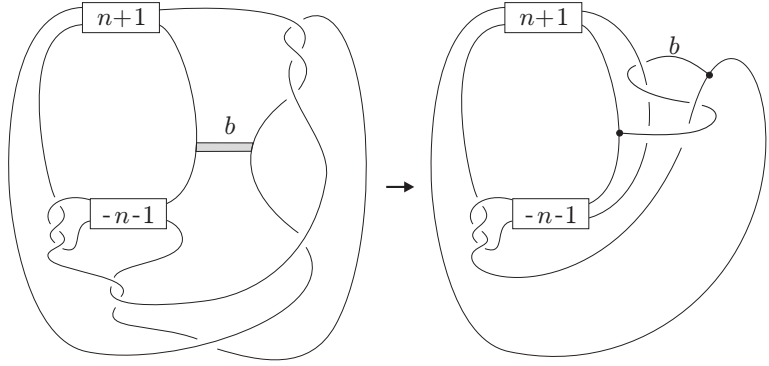}
\caption{$B_n(1/0)$ and the band $b$}
\label{fig:HHknot}
\end{figure}

\begin{figure}
\centering
\includegraphics[scale=0.65]{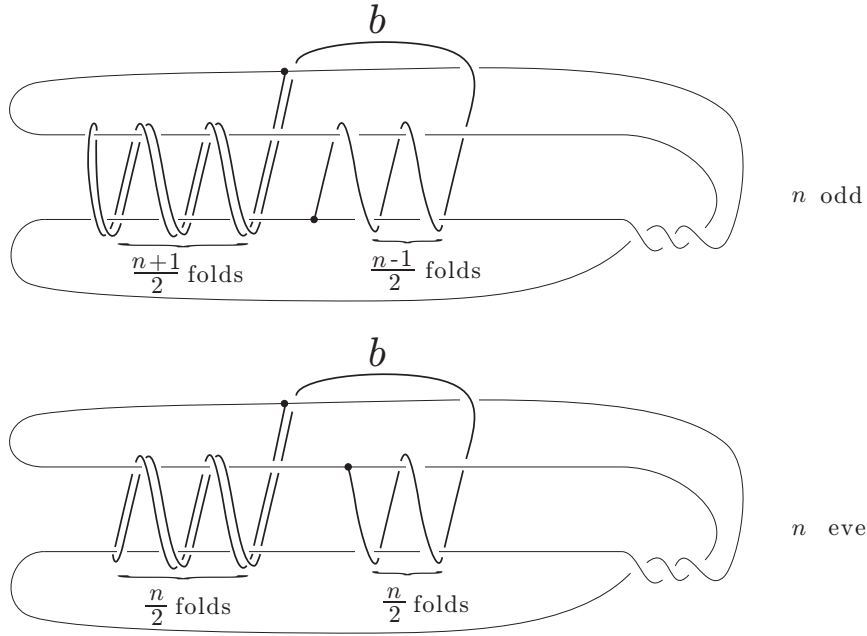}
\caption{$B_n(1/0)$ with $b$}
\label{fig:HHknot1}
\end{figure}

\begin{figure}
\centering
\includegraphics[scale=0.5]{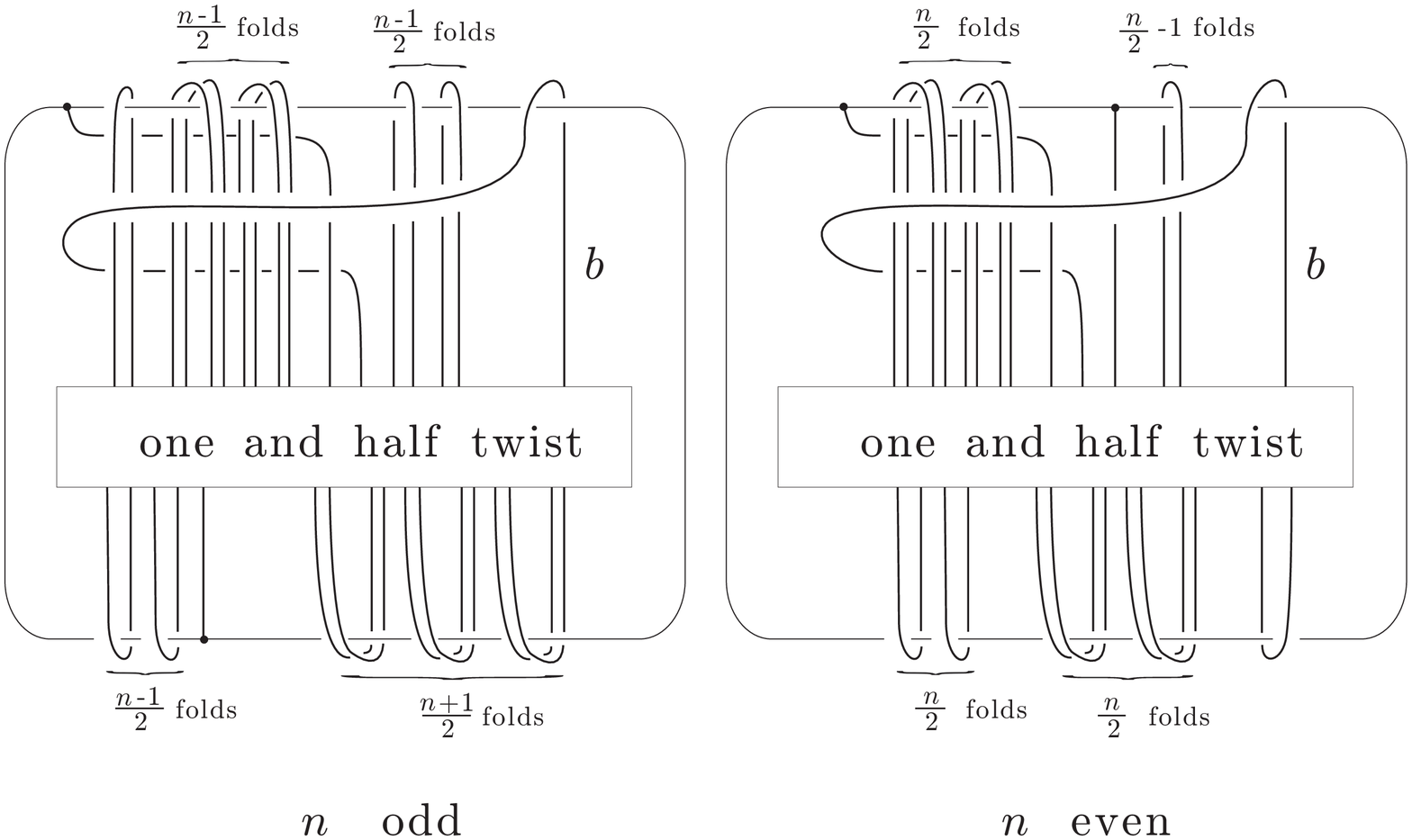}
\caption{The standard diagram of $B_n(1/0)$ with $b$}
\label{fig:HHknot2}
\end{figure}

\begin{proof}
Insert the $1/0$-tangle to $B_n$, and put a band $b$ as shown in Figure \ref{fig:HHknot} in order to keep track of framing.
Isotope the unknot $B_n(1/0)$ to a standard diagram as shown in Figure \ref{fig:HHknot2} (where the cases where $n=5$ and $n=4$ are drawn),
and take the double branched cover along it.
Then (the core of) the lift of $b$ gives $K_n$, and its framing corresponds to the $0$-tangle filling downstairs.
(In Figures \ref{fig:HHknot}, \ref{fig:HHknot1} and \ref{fig:HHknot2}, we draw $b$ in a line for simplicity during the deformation.)
From Figure \ref{fig:HHknot2},
we see that $K_n$ is the closure of a braid with $3n+2$ strings.
Moreover, there are $27n^2+41n+10$ positive crossings and $5n-1$ negative crossings.
After cancelling the negative crossings by positive crossings,
$K_n$ becomes the closure of a positive braid with $3n+2$ strings and $27n^2+36n+11$ crossings.
By \cite{St}, $K_n$ is fibered and Seifert's algorithm gives a fiber surface, whose genus is equal to the genus $g(K_n)$ of $K_n$.
Since $1-2g(K_n)=(3n+2)-(27n^2+36n+11)$, $g(K_n)=(27n^2+33n+10)/2$.

The framing of the lift of $b$ can be calculated to equal $m$.  This proves (2).
\end{proof}

Recall that $k^+(3n+1,3n+4)$ is hyperbolic for $n\ge 1$ by Lemma \ref{lem:k+1}.

\begin{proposition}\label{prop:HH}
For $n\ge 1$, let $K$ be the hyperbolic knot $K_n$ defined above, and
$K'$ be $k^+(3n+1,3n+4)$.
Let $m=27n^2+45n+21$.
Then $K$ and $K'$ are not equivalent, and $m$-surgery on $K$ and $K'$ yield homeomorphic lens spaces.
\end{proposition}

\begin{proof}
By Lemma \ref{lem:caseHH-Hlens}(1), $K$ has genus $(27n^2+33n+10)/2$, while
$K'$ has genus $(27n^2+33n+12)/2$ by Lemma \ref{lem:k+genus}.
Thus they are not equivalent.

Also, by Lemma \ref{lem:caseHH-Hlens}(2), $m$-surgery on $K$ yields $L(m,-9n^2-12n-5)=L(m,18n^2+33n+16)$.
As stated in Section \ref{sec:k+}, $m$-surgery on $K'$ yields $L(m,((3n+1)/(3n+4))^2)$.
Since
\begin{equation*}
\left(\frac{3n+1}{3n+4}\right)^2(18n^2+33n+16)\equiv 1 \pmod{m},
\end{equation*}
those lens spaces are homeomorphic.
\end{proof}

%%%%%%%%%%%%%%%%%%%%%%%%%%%%%%%%%%%%%%%%%%%%%%%%%%%%%%%%%%%%%%%%%%%%%%%%%%%%%%%%%%%%%
\section{Different classes}\label{sec:diff}

In this last section, we give the pairs of knots, each of which yields homeomorphic lens spaces by the same integral surgery,
and consists of knots belonging to different classes of hyperbolic, satellite, torus knots.

%%%%%%%%%
\subsection{Torus knot and satellite knot}

Let $C(a,b)$ be the $(2,2ab+1)$-cable of the torus knot of type $(a,b)$.

\begin{proposition}\label{prop:caseST}
For $n\ge 1$, let $K$ be the torus knot of type $(2n+1,4n+4)$, $K'=C(n+1,2n+1)$.
Let $m=8n^2+12n+5$.
Then $m$-surgery on $K$ and $K'$ yield homeomorphic lens spaces.
\end{proposition}

\begin{proof}
By \cite{Mo}, $m$-surgery on $K$ yield the lens space $L(m,(2n+1)^2)$.
Also, $m$-surgery on $K'$ yields $L(m,4(n+1)^2)$ by \cite{FS}.
Since $(2n+1)^2+4(n+1)^2=m$, these lens spaces are homeomorphic.
\end{proof}

%%%%%%%%%%
\subsection{Satellite knot and hyperbolic knot}

\begin{lemma}\label{lem:caseSH}
For $n\ge 0$, $4F_n^4+(-1)^nF_{n+2}^2=(4F_nF_{n+2}+(-1)^n)(F_{n+2}^2-4F_nF_{n+1})$.
\end{lemma}

\begin{proof}
First, 
$4F_n^4+(-1)^nF_{n+2}^2-(4F_nF_{n+2}+(-1)^n)(F_{n+2}^2-4F_nF_{n+1})= 4F_n(F_n^3-F_{n+2}^3+4F_nF_{n+1}F_{n+2}-(-1)^{n+1}F_{n+1})$.
From Cassini's identity,
\begin{eqnarray*}
\lefteqn{F_n^3-F_{n+2}^3+4F_nF_{n+1}F_{n+2}-(-1)^{n+1}F_{n+1}}\\
 &=& F_n^3-F_{n+2}^3+3F_nF_{n+1}F_{n+2}+F_{n+1}^3\\
 &=& F_n^3-(F_n+F_{n+1})^3+3F_nF_{n+1}F_{n+2}+F_{n+1}^3\\
 &=& -3F_nF_{n+1}(F_n+F_{n+1}-F_{n+2})=0.
\end{eqnarray*}
\end{proof}

Note that two Fibonacci numbers $F_n$ and $F_{n+2}$ are coprime, since
$\gcd(F_n,F_{n+2})=\gcd(F_n,F_{n+1})=1$.
By Lemma \ref{lem:k+2}, $k^+(F_{n+2},F_{n})$ is hyperbolic for $n\ge 3$.

\begin{proposition}\label{prop:caseSH}
For $n\ge 3$, let $K$ be the satellite knot $C(F_n,F_{n+2})$, $K'$ be the hyperbolic knot $k^+(F_{n+2},F_{n})$.
Let $m=4F_nF_{n+2}+(-1)^n$.
Then $m$-surgery on $K$ and $K'$ yield homeomorphic lens spaces.
\end{proposition}

\begin{proof}
By \cite{FS}, $m$-surgery on $K$ yields the lens space $L(m,4F_n^2)$.
From Lemma \ref{lem:F2}, $m$-surgery on $K'$ yields $L(m,(F_n/F_{n+2})^2)$.
Then
\[
4F_n^2 \left(\frac{F_n}{F_{n+2}}\right)^2\equiv (-1)^{n+1} \pmod{m}
\]
by Lemma \ref{lem:caseSH}(2), 
Thus the two lens spaces are homeomorphic.
\end{proof}

%%%%%%%%%
\subsection{Torus knot and hyperbolic knot}

\begin{figure}
\centering
\includegraphics[scale=0.5]{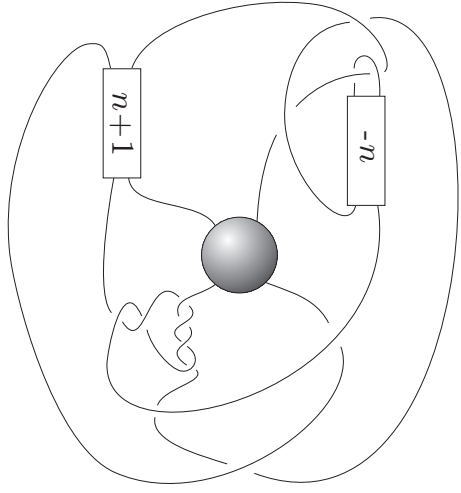}
\caption{The tangle $B_n$}
\label{fig:caseTH}
\end{figure}

For $n\ge 1$, let $B_n$ be the tangle as shown in Figure \ref{fig:caseTH},
where a vertical box denotes right-handed vertical half-twists.

Given $\alpha\in \mathbb{Q}\cup \{1/0\}$,
$B_n(\alpha)$ denotes the knot or link in $S^3$ obtained by inserting the rational tangle of slope $\alpha$
into the central puncture of $B_n$.
Also, $\widetilde{B}_n$ is the double branched cover of $S^3$ branched over $B_n(\alpha)$. 

\begin{lemma}\label{lem:caseTH-tangle}
\begin{enumerate}
\item $\widetilde{B}_n(1/0)=S^3$.
\item $\widetilde{B}_n(0)=L(18n^2+33n+15,18n+19)$.
\item $\widetilde{B}_n(-1)$ is a non-Seifert toroidal manifold $D^2(2,n+2)\cup D^2(4,2n+1)$.
\item $\widetilde{B}_n(1)$ is a non-Seifert toroidal manifold
$D^2(2,n)\cup D^2(5,2n+3)$ if $n\ge 2$, a Seifert fibered manifold of type $S^2(3,5,5)$ if $n=1$.
\end{enumerate}
\end{lemma}

\begin{figure}
\centering
\includegraphics[scale=0.4]{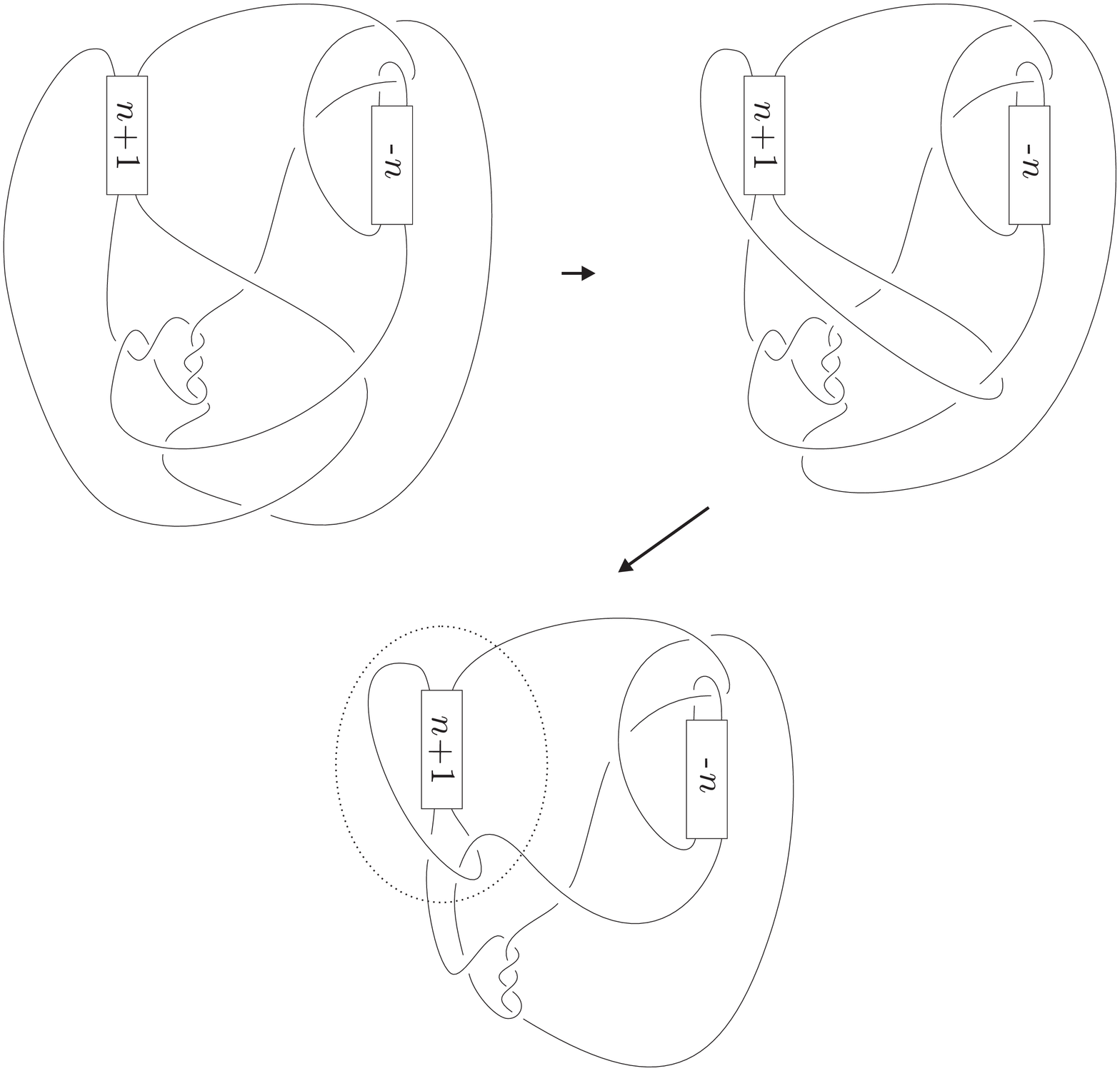}
\caption{$B_n(-1)$}
\label{fig:THknot1}
\end{figure}

\begin{figure}
\centering
\includegraphics[scale=0.4]{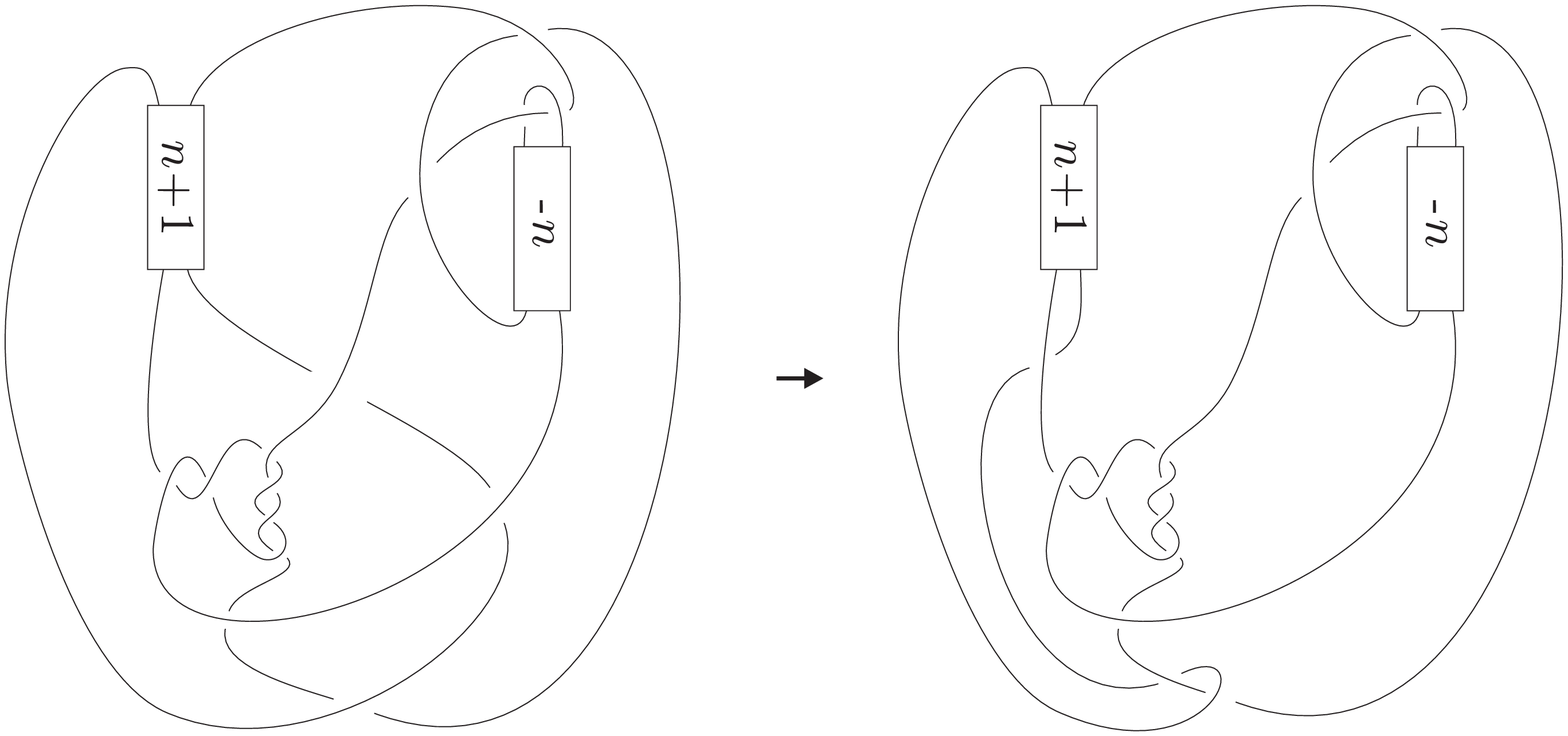}
\caption{$B_n(1)$}
\label{fig:THknot2}
\end{figure}

\begin{proof}
It is straightforward to see that $B(1/0)$ is the unknot and $B(0)$ is the $2$-bridge link corresponding to $(18n^2+33n+15)/(18n+19)$.
For $B(-1)$ and $B(1)$, see Figures \ref{fig:THknot1} and \ref{fig:THknot2}, respectively.
\end{proof}

By Lemma \ref{lem:caseTH-tangle}(1), the lift of $B_n$ in $\widetilde{B}_n(1/0)$ gives
the knot exterior of some knot $K_n$ in $S^3$, which is uniquely determined by Gordon-Luecke's theorem \cite{GL}.

\begin{lemma}
$K_n$ is hyperbolic.
\end{lemma}

\begin{proof}
This immediately follows from Lemmas \ref{lem:criterionH} and \ref{lem:caseTH-tangle}.
\end{proof}

\begin{lemma}\label{lem:caseTH-Hlens}
Let $m=18n^2+33n+15$.
Then $m$-surgery on $K_n$ yields the lens space $L(m,-18n-19)$. 
\end{lemma}

\begin{proof}
The argument is similar to the proof of Lemma \ref{lem:caseHH-Hlens}.
We omit it.
\end{proof}

\begin{proposition}\label{prop:caseTH}
For $n\ge 1$, let $K$ be the torus knot of type $(3n+2,6n+7)$,
$K'$ the knot $K_n$ defined above.
Let $m=18n^2+33n+15$.
Then $m$-surgery on $K$ and $K'$ yield homeomorphic lens spaces.
\end{proposition}

\begin{proof}
By \cite{Mo}, $m$-surgery on $K$ yields $L(m,9n^2+12n+4)$.
By Lemma \ref{lem:caseTH-Hlens}, $m$-surgery on $K'$ yields $L(m,18n+19)$.
Since $(9n^2+12n+4)(18n+19)\equiv 1 \pmod{m}$,
two lens spaces are homeomorphic.
\end{proof}

%%%%%%%%%%%%%%%%%%%%%%%%%%%
\begin{proof}[Proof of Theorem \textup{\ref{thm:main1}}]
This follows from Propositions \ref{prop:TT}, \ref{prop:HH}, \ref{prop:caseST}, \ref{prop:caseSH}
and \ref{prop:caseTH}.
\end{proof}

%%%%%%%%%%%%%%%%%%%%%%%%%%%%%%%%%%%%%%
We would like to thank Kazuhiro Kawasaki for computer experiments.

%%%%%%%%%%%%%%%%%%%%%%%%%%%%%%%%%%

\end{document}